\documentclass[a4paper,11pt]{article}

\usepackage{array}                   
\usepackage{theorem}                 
\usepackage{amsmath,amscd,amssymb}   
\usepackage{latexsym}
\usepackage{epsfig}

\theorembodyfont{\sl}

\newtheorem{lemma}{Lemma}[section]
\renewcommand{\thelemma}{\Roman{lemma}}
\newtheorem{proposition}[lemma]{Proposition}
\newtheorem{theorem}[lemma]{Theorem}
\newtheorem{corollary}[lemma]{Corollary}

\newcommand{\CC}{\mathbb C}

\newcommand{\HH}{\mathbb H}

\newcommand{\NN}{\mathbb N}

\newcommand{\QQ}{\mathbb Q}
\newcommand{\RR}{\mathbb R}

\newcommand{\ZZ}{\mathbb Z}
\newcommand{\cA}{\mathcal A}

\newcommand{\cH}{\mathcal H}

\newcommand{\cL}{\mathcal L}

\newcommand{\cO}{\mathcal O}
\newcommand{\cP}{\mathcal P}

\newcommand{\To}{\longrightarrow}
\newcommand{\Mapsto}{\mapstochar\longrightarrow}

\newcommand{\Sum}{\sum\limits}
\newcommand{\Prod}{\prod\limits}

\newcommand{\tens}{\otimes}

\newcommand{\contin}{\subseteq}

\renewcommand{\Bar}{\overline}

\newcommand{\cross}{\times}

\newcommand{\imic}{\cong}
\newcommand{\hf}{{\sfrac{1}{2}}}

\newcommand{\normal}{\vartriangleleft}
\newcommand{\Eins}{{\mathbf 1}}
\newcommand{\sqtimes}{{\boxtimes}}
\newcommand{\sfrac}[2]{{\textstyle{\frac{#1}{#2}}}}

\newcommand{\PSp}{\mathop{\mathrm {PSp}}\nolimits}

\newcommand{\PSL}{\mathop{\mathrm {PSL}}\nolimits}
\newcommand{\SL}{\mathop{\mathrm {SL}}\nolimits}
\newcommand{\Sp}{\mathop{\mathrm {Sp}}\nolimits}

\newcommand{\Pic}{\mathop{\mathrm {Pic}}\nolimits}

\newcommand{\diag}{\mathop{\mathrm {diag}}\nolimits}

\renewcommand{\Im}{\mathop{\mathrm {Im}}\nolimits}

\newcommand{\hcf}{\mathop{\mathrm {hcf}}\nolimits}
\newcommand{\red}{\mathop{\mathrm {red}}\nolimits}
\newcommand{\bda}{\mathbf a}
\newcommand{\bdb}{\mathbf b}
\newcommand{\bdc}{\mathbf c}

\newcommand{\bdr}{\mathbf r}

\newcommand{\bdv}{\mathbf v}
\newcommand{\bdw}{\mathbf w}
\newcommand{\bdx}{\mathbf x}
\newcommand{\bdy}{\mathbf y}

\newcommand{\gothS}{\mathfrak S}
\newcommand{\qedsymbol}{\mbox{$\Box$}}
\newcommand{\qed}{\unskip\nobreak\hfil\penalty50\hskip1em\hbox{}\nobreak
\hfill\qedsymbol\parfillskip=0pt\finalhyphendemerits=0}
\newenvironment{proof}{\begin{ProofwCaption}{Proof}}{\end{ProofwCaption}}
\newenvironment{ProofwCaption}[1]
  {\addvspace\theorempreskipamount \noindent{\it #1.}\rm}
  {\qed \par \addvspace\theorempostskipamount}

\newcommand{\Gtbl}{{\Gamma_t^{\mathrm{bil}}}}
\newcommand{\Gtnat}{{\Gamma_t^\natural}}
\newcommand{\Gtlev}{{\Gamma_t^{\mathrm{lev}}}}
\newcommand{\Gtblti}{{{\tilde\Gamma_t^{\mathrm{bil}}}}}
\newcommand{\Gtlevti}{{\tilde\Gamma_t^{\mathrm{lev}}}}
\newcommand{\Gtnatti}{{\tilde\Gamma_t^\natural}}

\newcommand{\Atbl}{{\cA_t^{\mathrm{bil}}}}
\newcommand{\Atcol}{{\cA_t^{\mathrm{col}}}}
\newcommand{\AtblV}{{\cA_t^{\mathrm{bil*}}}}

\newcommand{\Atlev}{{\cA_t^{\mathrm{lev}}}}
\newcommand{\AtlevV}{{\cA_t^{\mathrm{lev*}}}}
\newcommand{\ApblV}{{\cA_p^{\mathrm{bil*}}}}

\newcommand{\Aplev}{{\cA_p^{\mathrm{lev}}}}
\newcommand{\vz}{{\bdv_{(0,1)}}}
\begin{document}

\title{Abelian surfaces with odd bilevel structure}
\author{G.K. Sankaran}
\maketitle

Abelian surfaces with weak bilevel structure were introduced by
S.~Mukai in~\cite{Muk}. There is a coarse moduli space, denoted
$\Atbl$, for abelian surfaces of type $(1,t)$ with weak bilevel
structure. $\Atbl$ is a Siegel modular threefold, and can be
compactified in a standard way by Mumford's toroidal
method~\cite{AMRT}. We denote the toroidal compactification (in this
situation also known as the Igusa compactification) by $\AtblV$.
It is a projective variety over~$\CC$, and it is shown in~\cite{Muk} that
$\AtblV$ is rational for $t\leq 5$. In this paper we examine the
Kodaira dimension $\kappa(\AtblV)$ for larger~$t$. Our main result is
the following (Theorem~\ref{mainthm}).
\medskip

{\bf Theorem.} {\sl $\AtblV$ is of general type for $t$ odd and $t\geq 17$.}
\medskip

It follows from the theorem of L.~Borisov~\cite{Bor} that $\AtblV$ is
of general type for $t$ sufficiently large. If $t=p$ is prime, then it
follows from \cite{GH} and \cite{HS} that $\ApblV$ is of general type
for $p\geq 37$. Our result provides an effective bound in the general
case and a better bound in the case $t=p$. As far as we know, all
previous explicit general type results (for instance
\cite{GH,HS,O'G,GS,San}) have been for the cases $t=p$ or $t=p^2$ only.

It is for brevity that we assume $t$ is odd. If $t$ is even the
combinatorial details are more complicated, especially when
$t\equiv 2$ mod~$4$, but the method is still applicable. In fact the method
is essentially that of~\cite{HS}, with some modifications.

{\it Acknowledgement\/}. Part of this work resulted from
conversations with my research student Alfio Marini.

\section{Background}
\renewcommand{\thelemma}{\Roman{section}.\arabic{lemma}}

If $A$ is an abelian surface with a polarisation $H$ of type $(1,t)$,
$t>1$, then a {\em canonical level structure}, or simply {\em level
structure}, is a symplectic isomorphism
\begin{equation*}\label{definelevelstr}
\alpha:\ZZ_t^2\To K(H)=\left\{\bdx\in A\mid t_\bdx^*\cL\imic\cL\mbox{
if }c_1(\cL)=H\right\}.
\end{equation*}
The moduli space $\Atlev$ of abelian surfaces with a canonical level
structure has been studied in detail in~\cite{HKW2}, chiefly in the
case~$t=p$.

A {\em colevel structure} on~$A$ is a level structure on the dual
abelian surface $\hat A$: note that $H$ induces a polarisation $\hat
H$ on $\hat A$, also of type~$(1,t)$. Alternatively, a colevel
structure may be thought of as a symplectic isomorphism
\begin{equation*}\label{definecolevelstr}
\beta:\ZZ_t^2\To A[t]/K(H)
\end{equation*}
where $A[t]$ is the group of all $t$-torsion points of~$A$. Obviously
the moduli space $\Atcol$ of abelian surfaces of type $(1,t)$ with a
colevel structure is isomorphic to $\Atlev$, and each of them has a
forgetful morphism $\psi^{\mathrm{lev}}$, $\psi^{\mathrm{col}}$ to the
moduli space $\cA_t$ of abelian surfaces of type~$(1,t)$. We define
\begin{equation*}\label{defineAtlev}
\Atbl=\Atlev\cross_{\cA_t}\Atcol.
\end{equation*}
The forgetful map $\psi^{\mathrm{lev}}:\Atlev\to\cA_t$ is the quotient
map under the action of $\SL(2,\ZZ_t)$ given by
\begin{equation*}\label{defineSL2action}
\gamma:\left[(A,H,\alpha)\right]\mapsto\left[(A,H,\alpha\gamma)\right]
\end{equation*}
where $\gamma\in\SL(2,\ZZ_t)$ is viewed as a symplectic automorphism
of~$\ZZ_t^2$. The action is not effective, because $(A,H,\alpha)$ is
isomorphic to $(A,H,-\alpha)$ via the isomorphism $\bdx\mapsto -\bdx$;
so $-\Eins_2\in\SL(2,\ZZ_t)$ acts trivially. Thus $\psi^{\mathrm{lev}}$ is a
Galois morphism with Galois group~$\PSL(2,\ZZ_t)=\SL(2,\ZZ_t)/\pm\Eins_2$.

A point of $\Atbl$ thus corresponds to an equivalence class
$\left[(A,H,\alpha,\beta)\right]$, where $(A,H)$ is a polarised
abelian surface of type~$(1,t)$, $\alpha$ and $\beta$ are level and
colevel structures, and $(A,H,\alpha,\beta)$ is equivalent to
$(A',H',\alpha',\beta')$ if there is an isomorphism $\rho:A\to A'$
such that $\rho^* H'=H$, $\rho\alpha=\alpha'$ and
$\hat\rho^{-1}\beta=\beta'$. In particular, for general $A$, we have
$(A,H,\alpha,\beta)\imic(A,H,-\alpha,-\beta)$ but
$(A,H,\alpha,\beta)\not\imic(A,H,-\alpha,\beta)$. Another way to
express this is to say that the wreath product
$\ZZ_2\wr\PSL(2,\ZZ_t)$, acts on $\Atbl$ with quotient~$\cA_t$.
\begin{theorem}\label{proveGtbl}
 {\rm (Mukai~\cite{Muk})} $\Atbl$ is the quotient of the Siegel upper
half-plane $\HH_2$ by the group
\begin{equation*}
\Gtbl=\Gtnat\cup\zeta\Gtnat
\end{equation*}
where
\begin{equation*}\label{defineGtbl}
\Gtnat=\left\{\gamma\in\Sp(4,\ZZ)\mid \gamma-\Eins_4\in
\begin{pmatrix}t\ZZ&\ast&t\ZZ&t\ZZ\\
                t\ZZ&t\ZZ&t\ZZ&t^2\ZZ\\
                t\ZZ&\ast&t\ZZ&t\ZZ\\
                \ast&\ast&\ast&t\ZZ
\end{pmatrix}
\right\}
\end{equation*}
and $\zeta=\diag(1,-1,1,-1)$, acting by fractional linear transformations.
\end{theorem}

Thus $\Gtbl$ should be thought of as a subgroup of the paramodular
group
\begin{equation*}\label{defineparamodular}
\Gamma_t=\left\{\gamma\in\Sp(4,\QQ)\mid \gamma-\Eins_4\in
\begin{pmatrix}\ast&\ast&\ast&t\ZZ\\
                t\ZZ&\ast&t\ZZ&t\ZZ\\
                \ast&\ast&\ast&t\ZZ\\
                \ast&\frac{1}
{t}\ZZ&\ast&\ast
\end{pmatrix}
\right\}.
\end{equation*}
(The paramodular group is the group denoted
$\Gamma_{1,t}^\circ$ in~\cite{HKW2} and~\cite{FS}.)

For some purposes it is more convenient to work with the conjugate
$\Gtblti=R_t\Gtbl R_{t^{-1}}$ of $\Gtbl$ by $R_t=\diag(1,1,1,t)$,
and with the corresponding
conjugates $\Gtnatti$, $\Gtlevti$ etcetera. These groups have the
advantage that they are subgroups of
$\Sp(4,\ZZ)$ rather than $\Sp(4,\QQ)$, and are defined by congruences
mod~$t$, not mod~$t^2$, but their action on $\HH_2$ is not the usual one by
fractional linear transformations.

If $E_i$ are elliptic curves and $(A,H)=\left(E_1\cross E_2,
c_1\big(\cO_{E_1}(1)\sqtimes\cO_{E_2}(t)\big)\right)$, we say that
$(A,H)$ is a product surface. In this case $K(H)=\{0_{E_1}\}\cross
E_2[t]$, so a level structure on $A$ may be thought of as a full
level-$t$ structure on~$E_2$. The automorphism $(\bdx,\bdy)\mapsto
(\bdx,-\bdy)$ of $A=E_1\cross E_2$ induces an isomorphism
$(A,H,\alpha,\beta)\to(A,H,-\alpha,\beta)$ in this case, so a product
surface with a weak bilevel structure still has an extra automorphism. The
corresponding locus in the moduli space arises from the fixed locus of
$\zeta$ in $\HH_2$, and will be of great importance in this paper.

The geometry of $\AtblV$ shows many similarities with that of
$\AtlevV$, which was studied (in the case of $t$ an odd prime) in the
book~\cite{HKW2}. In many cases where the proofs of intermediate
results are very similar to those of corresponding results
in~\cite{HKW2} we omit the details and simply indicate the appropriate
reference.

\section{Modular groups and modular forms}

We first collect some facts about congruence subgroups in $\SL(2,\ZZ)$
and some related combinatorial information. For $r\in\NN$ we denote by
$\Gamma_1(r)$ the principal congruence subgroup of $\SL(2,\ZZ)$. We
denote the modular curve $\Gamma_1(r)\backslash\HH$ by $X^\circ(r)$,
and the compactification obtained by adding the cusps by~$X(r)$.

For $m, r\in\NN$, define
\begin{equation*}
\Phi_m(r)=\{\bda\in \ZZ_r^m\mid
\bda\mbox{ is not a multiple of a zerodivisor in }\ZZ_r\},
\end{equation*}
that is, $\bda\in\Phi_m(r)$ if and only if $\bda=z\bda'$ implies
$z\in\ZZ_r^*$; and put $\phi_m(r)=\#\Phi_m(r)$. We also put
$\Bar\Phi_m(r)=\Phi_m(r)/\pm 1$.
\begin{lemma}\label{calculatePhimt}
If the primes dividing $r$ are $p_1<p_2<\ldots<p_n$ then 
\begin{equation*}
\phi_m(r)=\sum_{i=0}^n(-1)^i\sum_{p_{j_1},\dots,p_{j_i}}
\bigg(r\prod_{k=1}^i p_{j_k}^{-1}\bigg)^m=r^m\prod_{p|r}(1-p^{-m}).
\end{equation*}
\end{lemma}
\begin{proof}
We first prove that $\phi_m(r)$ is a multiplicative function. First we
suppose that $r=pq$, with $\gcd(p,q)=1$. It is easy to see that
$\bda\in \Phi_m(r)$ if and only if $\bda_p\in\Phi_m(p)$ and 
$\bda_q\in\Phi_m(q)$, where $\bda_p$ denotes the reduction of $\bda$
mod~$p$.

We divide $\ZZ_r^m$ into residue classes mod~$p$: that is, we write
$\ZZ_r^m$ as the disjoint union of subsets $S_\bdc$ for
$\bdc\in\ZZ_p^m$, where $S_\bdc=\{\bda\mid \bda_p=\bdc\}$. There are 
$\phi_m(p)$ subsets~$S_\bdc$ such that $\bdr\in\Phi_m(p)$.

The reduction mod~$q$ map $S_\bdc\to \ZZ_q^m$ is bijective, since it
is the inverse of the injective map $\bdb\mapsto
\bdc+p\bdb\in\ZZ_r^m$. Hence in each of the $\phi_m(p)$ subsets
$S_\bdc$, $\bdc\in\Phi_m(p)$ there are $\phi_m(q)$ elements whose
reduction mod~$q$ belongs to $\Phi_m(q)$. It follows that
$\phi_m(r)=\phi_m(p)\phi_m(q)$.

Finally, we check that if $r=p^k$, $p$~prime, then
$\phi_m(r)=r^m(1-p^{-m})$. If $\bda\not\in\Phi_m(r)$, then
$\bda=p\bda'$ for a unique $\bda'\in\ZZ_{r/p}^m$, so there are
$(p^{k-1})^m$ such elements~$\bda$.
\end{proof}
Note that $\phi_1$ is the Euler $\phi$ function,
and $\Phi_1(r)$ is the set of non-zerodivisors of $\ZZ_r$. 
\begin{corollary}\label{orderofSL2}
The order of $\SL(2,\ZZ_t)$ is given by
\begin{equation*}\label{orderofSLformula}
|\SL(2,\ZZ_t)|=t\phi_2(t)=t^3\prod_{p|t}(1-p^{-2}).
\end{equation*}
\end{corollary}

\begin{proof}
(See also~\cite[\S 1.6]{Shi}.) If $A\in \SL(2,\ZZ_t)$, then
$A_1=(a_{11},a_{12})\in \Phi_2(t)$. So by Euclid's algorithm we can
find $A'_2=(a'_{21},a'_{22})$ such that
$\det\begin{pmatrix}A_1\\ A'_2\end{pmatrix}=\gcd(a_{11},a_{12})=r$.
Replacing $A'_2$ by $A_2=r^{-1} A'_2$, we get a matrix $A$ with
$\det(A)=1$.  Furthermore, if
$B_j=\begin{pmatrix}A_1\\A_2+jA_1\end{pmatrix}$,
$j=0,\dots,t-1$,
then $\det(B_j)=\det(A)=1$, and $B_j\neq B_{j'}$ if $j\neq j'$.  So
$|\SL(2,\ZZ_t)|=t\phi_2(t)$.
\end{proof}

For $r>2$, put $\mu(r)=\left[\PSL(2,\ZZ):\Gamma_1(r)\right]$. By
Corollary~\ref{orderofSL2} we have
\begin{equation*}\label{calculatemu}
\mu(r)=r^3\prod_{p|r}(1-p^{-2}).
\end{equation*}
We need the following well-known lemma.
\begin{lemma}\label{cuspsandgenusofX(t)}
If $r>2$ then $X(r)$ has 
\begin{equation*}
\nu(r)=\mu(r)/r=r^2\prod_{p|r}(1-p^{-2})
\end{equation*}
cusps and is a smooth complete curve of genus
$g=1+\frac{\mu(r)}{12}-\frac{\nu(r)}{2}$.
\end{lemma}
\begin{proof}
See~\cite[pp.~23--24]{Shi}.
\end{proof}
We denote $\mu(t)$ by $\mu$ and $\nu(t)$ by~$\nu$. Note that
$\phi_2(1)=\nu(1)=1$ and $\phi_2(r)=2\nu(r)$ for~$r>2$.

Now we turn to subgroups of $\Sp(4,\QQ)$ and to modular forms. Denote
by $\gothS_n^*(\Gamma)$ the space of weight~$n$ cusp forms
for~$\Gamma\subseteq \Sp(4,\QQ)$. We need the groups
$\bar\Gamma(1)=\PSp(4,\ZZ)$ and, for $\ell\in\NN$,
\begin{equation*}\label{definefulllevel}
\Gamma(\ell)=\big\{\gamma\in\Sp(4,\ZZ)\mid
\bar\gamma=\Eins_4\in\Sp(4,\ZZ_\ell)\big\}.
\end{equation*}
If $t^2|\ell$ then $\Gamma(\ell)\vartriangleleft\Gtbl$, because
$\Gamma(\ell)\contin\Gtbl$ and $\Gamma(\ell)$ is normal in
$\Gamma(1)=\Sp(4,\ZZ)$.

By a previous calculation~\cite{Tai} we know that
\begin{equation*}\label{estimatecuspforms}
\dim\gothS_n^*\big(\Gamma(\ell)\big) =
{\frac{n^3}{8640}}\left[\bar\Gamma(1):\Gamma(\ell)\right]+O(n^2)
\end{equation*}
(as long as $\ell>2$ we can consider $\Gamma(\ell)$ as a subgroup of
$\PSp(4,\ZZ)$ rather than $\Sp(4,\ZZ)$). A standard application of the
Atiyah--Bott fixed-point theorem (see~\cite{Hir}, or in this
context~\cite{HS}) gives
\begin{equation*}\label{atiyahbottgeneral}
\dim\gothS_n^*\big(\Gtbl\big)={\frac{a}{\left[\Gtbl:\Gamma(\ell)\right]}}
\dim\gothS_n^*\big(\Gamma(\ell)\big) +O(n^2)
\end{equation*}
where $a$ is the number of elements $\gamma\in\Gtbl$ whose fixed locus
in $\HH_2$ has dimension~$3$. Thus $a$ is the number of elements of
$\Gtbl$ that act trivially on~$\HH_2$. In $\Sp(4,\ZZ)$ there are two
such elements, $\pm\Eins_4$, but if $t>2$ then
$-\Eins_4\not\in\Gtbl$. So $a=1$, and hence
\begin{eqnarray}\label{cforms}
\dim\gothS_n^*\big(\Gtbl\big)&=&{\frac{1}{\left[\Gtbl:\Gamma(\ell)\right]}}
\dim\gothS_n^*\big(\Gamma(\ell)\big) +O(n^2) \nonumber\\
      &=&{\frac{n^3}{8640}}{\frac{\left[\bar\Gamma(1):\Gamma(\ell)\right]}
      {\left[\Gtbl:\Gamma(\ell)\right]}}+O(n^2)\nonumber\\
      &=&{\frac{n^3}{8640}}\left[\bar\Gamma(1):\Gtbl\right]+O(n^2).
\end{eqnarray}
The number $\left[\bar\Gamma(1):\Gtbl\right]$ is equal to the degree
of the map $\Atbl\to\cA_1$ (actually there are two such maps
of the same degree), where $\cA_1$ is the moduli space of principally
polarized abelian surfaces. Now
\begin{eqnarray*}
\left[\bar\Gamma(1):\Gtbl\right]&=&\hf\left[\bar\Gamma(1):\Gtnat\right]\\
 &=&\hf\left[\bar\Gamma(1):\Gtlev\right]\left[\Gtlev:\Gtnat\right]. 
\end{eqnarray*}
We can see directly that $\Gtlev\supset\Gtnat$ since
\begin{equation*}\label{defineGtlev}
\Gtlev=\left\{\gamma\in\Sp(4,\ZZ)\mid \gamma-\Eins_4\in
\begin{pmatrix}
         \ast&\ast&\ast&t\ZZ\\
         t\ZZ&t\ZZ&t\ZZ&t^2\ZZ\\
         \ast&\ast&\ast&t\ZZ\\
         \ast&\ast&\ast&t\ZZ
\end{pmatrix}
\right\}.
\end{equation*}

\begin{lemma}\label{GtlevontoSL2t}
The map 
\begin{equation*}
\varphi:\Gtlev\To \SL(2,\ZZ_t),\;A \mapsto
\begin{pmatrix}
a_{11}&a_{13}\\
a_{31}&a_{33}
\end{pmatrix}
\end{equation*}
is a surjective group homomorphism, and the kernel is $\Gtnat$. 
\end{lemma}
\begin{proof}
The surjectivity follows from the well-known fact that the redution
mod~$t$ map $\red_t:\SL(2,\ZZ)\to\SL(2,\ZZ_t)$ is surjective, and the
rest is obvious.
\end{proof}

\begin{lemma}\label{indexofGtlev}
For $t>2$, the index $[\bar\Gamma(1): \Gtlev]$ is equal to
$t\phi_4(t)/2$. 
\end{lemma}
\begin{proof}
The proof is almost the same as proof of \cite[Lemma 0.5]{HW}. In
place of the chain of groups $\Gamma_{1,p} <\:_0\Gamma_{1,p}
<\Gamma'=\Gamma(1)$, we use the chain $\Gtlev< \:_0\Gamma_{1,t}<
\Gamma(1)$. Furthermore, we use the set $\Phi_4(t)$ where
$\SL(4,\ZZ_t)$ acts. Note that $\SL(4,\ZZ)$ still acts transitively on
$\Phi_4(t)$, via
\begin{equation*}
\begin{pmatrix}
b_{11}&0&b_{12}&0\\
0&1&0&0\\
b_{21}&0&b_{22}&0\\
0&0&0&1
\end{pmatrix}
\mbox{ and }
\begin{pmatrix}
B&0\\
0&^tB^{-1}
\end{pmatrix},
\end{equation*}
for $B\in \SL(2,\ZZ)$.

Following the same steps as in~\cite{HW}, and substituting $\phi_m(t)$
for $p^m-1=\phi_m(p)$, we then find that
$[\:{}_0\Gamma_{1,t}:\Gtlev]=t\phi_1(t)$ and
$[\:{}_0\Gamma_{1,t}:\Gamma(1)|=\phi_4(t)/\phi_1(t)$,
so $[\overline{\Gamma}(1):\Gtlev]=t\phi_4(t)/2$.
\end{proof}

\begin{theorem}\label{cuspforms}
The number of cusp forms of weight~$n$ for $\Gtbl$ (for $t>2$) is given by
\begin{eqnarray*}
\dim \gothS_n^*(\Gtbl)&=&{\frac{n^3}{34560}}t^2\phi_2(t)\phi_4(t)
\label{cuspformsformula}\\ 
&=&{\frac{n^3}{34560}}t^8\prod_{p|t}(1-p^{-2})(1-p^{-4}).
\label{cuspformsproduct}
\end{eqnarray*} 
\end{theorem}
                                                         
\begin{proof} Immediate from equation~\eqref{cforms},
Corollary~\ref{orderofSL2} and Lemma~\ref{indexofGtlev}.
\end{proof}

\section{Torsion in the modular group}

We know that $\Gtbl\subset\Sp(4,\ZZ)$, and the conjugacy classes of
torsion elements in $\Sp(4,\ZZ)$ are known
(\cite{Gott,Ueno}). See~\cite{HKW1} for a summary of the relevant
information.

If $\gamma\in\Gtnat$ then the reduction mod~$t$ of $\gamma$ is 
\begin{equation*}\label{reduceGtnat}
\bar\gamma=\begin{pmatrix}
       1&\ast&0&0\\
       0&1&0&0\\
       0&\ast&1&0\\
       \ast&\ast&\ast&1
\end{pmatrix}
\in\Sp(4,\ZZ_t),
\end{equation*}
so the characteristic polynomial $\chi(\bar\gamma)$ is
$(1-x)^4\in\ZZ_t[x]$. On the other hand, if
$\gamma\in\zeta\Gtnat$ then
\begin{equation*}\label{reducezetaGtnat}
\bar\gamma=\zeta\begin{pmatrix}
       1&\ast&0&0\\
       0&1&0&0\\
       0&\ast&1&0\\
       \ast&\ast&\ast&1
\end{pmatrix}
=\begin{pmatrix}
          1&\ast&0&0\\
          0&-1&0&0\\
          0&\ast&1&0\\
          \ast&\ast&\ast&-1
\end{pmatrix}
\in\Sp(4,\ZZ_t),
\end{equation*}
so $\chi(\bar\gamma)=(1-x)^2(1+x)^2\in\ZZ_t[x]$.

The only classes in the list in \cite{Ueno}, up to conjugacy, where
the characteristic polynomials have this reduction mod~$t$ ($t>2$) are
I(1), where $\chi(\gamma)=(1-x)^4$, II(1)a and II(1)b. Class
I(1) consists of the identity; class II(1)a includes $\zeta$ so this
just gives us the conjugacy class of $\zeta$. Class II(2)b is the
$\Sp(4,\ZZ)$-conjugacy class of $\xi$, where
\begin{equation*}\label{xi}
\xi=\begin{pmatrix}
       1& 1& 0&0\\
       0&-1& 0&0\\
       0& 0& 1&0\\
       0& 0& 1&-1
\end{pmatrix} \in\Gtbl.
\end{equation*}

\begin{proposition}\label{torsion}
Every nontrivial element of finite order in $\Gtbl$ (for $t>2$) has
order~$2$, and is conjugate to $\zeta$ or to $\xi$ in $\Gtbl$ if $t$
is odd.
\end{proposition}

\begin{proof} It follows from the list in~\cite{Ueno} that the only torsion for
$t>2$ is $2$-torsion (this is still true if $t$ is even). The $2$-torsion
of the group $\Gtlev$ was studied by Brasch~\cite{Br}. There are five types
but only two of them occur for odd~$t$. The representatives for these
conjugacy classes given in~\cite{Br} are (up to sign) $\zeta$ and~$\xi$; so
the assertion of the theorem is that the $\Gtbl$-conjugacy classes of
$\zeta$ and $\xi$ coincide with the intersections of their $\Gtlev$-conjugacy
classes with~$\Gtbl$. This is checked in~\cite[Proposition~3.2]{SS} for
the case $t=6$ (the relevant cases are called $\zeta_0$ and $\zeta_3$ there),
but the proof is valid for all~$t>2$.
\end{proof}

We put
\begin{equation}\label{definecH1}
\cH_1=\left.\left\{\begin{pmatrix}\tau_1&0\\ 0&\tau_3\end{pmatrix}\,\right|\,
 \Im\tau_1>0, \;\Im\tau_3>0\right\}\subset\HH_2
\end{equation}
and
\begin{equation}\label{definecH2}
\cH_2=\left.\left\{\begin{pmatrix}\tau_1&\tau_2\\ \tau_2&\tau_3\end{pmatrix}
\,\right|\, 2\tau_2+\tau_3=0\right\}\subset\HH_2.
\end{equation}
These are the fixed loci of $\zeta$ and $\xi$ respectively. We denote by
$H^\circ_1$ and $H^\circ_2$ the images of $\cH_1$ and $\cH_2$ in 
$\Atbl$, and by $H_1$ and $H_2$ their respective closures in~$\AtblV$.
\begin{lemma}\label{irreducibleH}
$H^\circ_i$ is irreducible for $i=1,2$.
\end{lemma}                
\begin{proof} 
This follows at once from Proposition~\ref{torsion} together with
equations~\eqref{definecH1} and~\eqref{definecH2}.
\end{proof}

The abelian surfaces corresponding to points in $H^\circ_1$ and
$H^\circ_2$ are, respectively, product surfaces and bielliptic abelian
surfaces, as described in~\cite{HW} for the case $t$ prime.

We define the subgroup $\Gamma(2t,2t)$ of $\Gamma(t)\cross\Gamma(t)$ by
\begin{equation*}
\Gamma(2t,2t)=\{(M,N)\in \Gamma(t)\cross\Gamma(t)\mid M\equiv {}^\top\!
N^{-1} \mod 2\}
\end{equation*}
\begin{lemma}
$H^\circ_1$ is isomorphic to $X^\circ(t)\cross X^\circ(t)$, and $H^\circ_2$ is
isomorphic to $\Gamma(2t,2t)\backslash\HH\cross\HH$.
\end{lemma}
\begin{proof}
Identical to the proofs of the corresponding results \cite[Lemma
I.5.43]{HKW2} and \cite[Lemma I.5.45]{HKW2}. The level-$t$ structure
now occurs in both factors, whereas in~\cite{HKW2} there is level-$1$
structure in the first factor and level-$p$ structure in the
second. In~\cite{HKW2} the level $p$ is assumed to be an odd prime but this
fact is not used at that stage: $p$ odd suffices, so we may replace
$p$ by $t$. Thereafter one simply replaces all the groups with their
intersection with $\Gtbl$, which imposes a level-$t$ structure in the
first factor and causes it to behave exactly like the second factor.
\end{proof}
\begin{lemma}
$H^\circ_1$ and $H^\circ_2$ are disjoint.
\end{lemma}
\begin{proof}
The stabiliser of any point of $\HH_2$ in $\Gtbl$ is cyclic (of
order~$2$), since $\Gtnat$ is torsion-free and therefore has no fixed
points. A point of $\cH_1\cap\cH_2$ would be the image of a point of
$\HH_2$ stabilised by the subgroup generated by $\zeta$ and $\xi$,
which is not cyclic.
\end{proof}

\section{Boundary divisors}

We begin by counting the boundary divisors. These correspond to
$\Gtblti$-orbits of lines in $\QQ^{\,4}$: we identify a line by its
unique (up to sign) primitive generator
$\bdv=(v_1,v_2,v_3,v_4)\in\ZZ^4$ with $\hcf(v_1,v_2,v_3,v_4) = 1$. We
denote the reduction of $\bdv$ mod~$t$ by $\bar\bdv=(\bar v_1,\bar
v_2,\bar v_3,\bar v_4)\in\ZZ_t^4$. To fix things we shall say,
arbitrarily, that $\bdv$ is positive if the first non-zero entry $\bar
v_i$ of $\bar\bdv$ satisfies $\bar v_i\in\{1,\ldots,(t-1)/2\}$
(remember that we have assumed that $t$ is odd). Then each line has a
unique positive primitive generator.

If $\bdv=(v_1,v_2,v_3,v_4)\in\ZZ^4$, we define the $t$-divisor to be
$r=\hcf(t,v_1,v_3)$.

\begin{proposition}\label{bdydivisors}
The lines $\QQ\,\bdv$ and $\QQ\,\bdw$ spanned by positive primitive vectors
$\bdv,\bdw\in\ZZ^4$ are in the
same $\Gtblti$-orbit if and only if $(\bar v_1,\bar v_3)=(\bar
w_1,\bar w_3)$ (in particular $\bdv$ and $\bdw$ have the same
$t$-divisor, $r$), and $(v_2, v_4)\equiv\pm(w_2,w_4)$ mod~$r$.
\end{proposition}
\begin{proof} 
Note that if $\Gamma(t)$ is the principal congruence subgroup of
level~$t$ in $\Sp(4,\ZZ)$ then $\Gamma(t)\normal\Gtnatti$ and the
quotient is
\begin{equation*}
\Gtnatti(t)=\left\{
\begin{pmatrix}
1&k&0&k'\\ 0&1&0&0\\ 0&l&1&l'\\ 0&0&0&1
\end{pmatrix}
\in\Sp(4,\ZZ_t)\right\}
\imic \ZZ_t^4.
\end{equation*}
We claim that two primitive vectors $\bdv$ and $\bdw$ are equivalent
modulo $\Gamma(t)$ if and only if $\bar v=\bar w$. It is obvious that
$\Gamma(t)$ preserves the residue classes mod~$t$. Conversely, suppose
that $\bar v=\bar w$. Then we can find $\gamma\in\Sp(4,\ZZ)$ such that
$\gamma\bdv=(1,0,0,0)$ (the corresponding geometric fact is that the
moduli space $\cA_2$ of principally polarised abelian surfaces has
only one rank~$1$ cusp). Since $\Gamma(t)\normal\Sp(4,\ZZ)$ this means
that in order to prove the claim we may assume $\bdv=(1,0,0,0)$. Then
we proceed exactly as in the proof of \cite[Lemma 3.3]{FS}, taking
$p=1$ and $q=t$ (the assumptions that $p$ and $q$ are prime are not
used at that point).

The group $\Gtnatti(t)$ acts on the set $(\ZZ_t^4)^\cross$ of non-zero
elements of $\ZZ_t^4$ by $\bar v_2\mapsto \bar v_2+k\bar v_1 + l\bar
v_3$ and $\bar v_4\mapsto \bar v_4+k'\bar v_1 + l'\bar v_3$: so
$\bar\bdv$ is equivalent to $\bar\bdw$ if and only if $(\bar v_1,\bar
v_3)=(\bar w_1,\bar w_3)$, so they have the same $t$-divisor, and
$\bar v_2\in\bar w_2+\ZZ_t r$ and $\bar v_4\in\bar w_4+\ZZ_t r$. These
are therefore the conditions for primitive vectors $\bdv$ and $\bdw$
to be equivalent under~$\Gtnatti$. For equivalence under $\Gtblti$,
we get the extra element~$\zeta$ which makes $(v_1,v_2,v_3,v_4)$
equivalent to $(v_1,-v_2,v_3,-v_4)$. Since we are interested in orbits
of lines, not primitive generators, we may restrict ourselves to
positive generators~$\bdv$.
\end{proof}

The irreducible components of the boundary divisor of $\AtblV$
correspond to the $\Gtbl$-orbits (or equivalently to
$\Gtblti$-orbits) of lines in $\QQ^{\,4}$. We denote the boundary
component corresponding to $\QQ\,\bdv$ by $D_\bdv$. We shall be
chiefly interested in the cases $r=t$ and $r=1$. We refer to these as
the standard components. They are represented by vectors $(0,a,0,b)$
and $(a,0,b,0)$ respectively, in both cases with $\hcf(a,b)=1$, $0\le
a \le (t-1)/2$ and $0\le b <t$. Note that there are $\nu$ of each of these.

\begin{corollary}\label{numberofbdycpts}
If $t$ is odd then the number of irreducible boundary divisors of $\AtblV$
with $t$-divisor $r$ is $\#\Bar\Phi_2(h)\#\Bar\Phi_2(r)$, where $h=t/r$. For
$r\neq 1$, $t$, this is equal to $\frac{1}{4}\phi_2(h)\phi_2(r)$.
\end{corollary}
\begin{proof}
See above for the standard cases. In general, the $\Gtnat$-orbit of a
primitive vector $\bdv$ is determined by the classes of 
$(v_1/r,v_3/r)$ in $\Phi_2(h)$ and of $(\bar v_2,\bar v_4)\in \Phi_2(r)$.
The extra element $\zeta$ and the freedom to multiply $\bdv$
by $-1\in\QQ$ allow us to multiply either of these classes by~$-1$ and
the choices therefore lie in $\bar\Phi_2(h)$ and $\Bar\Phi_2(r)$.
\end{proof}

\section{Jacobi forms}

In this section we shall describe the behaviour of a modular form
$F\in\gothS_{3n}^*(\Gtbl)$ near a boundary divisor $D_\bdv$. The
standard boundary divisors are best treated separately, since it is in those
cases only that the torsion plays a role: on the other hand,
the standard boundary divisors occur for all~$t$ and their behaviour
is not much dependent on the factorisation of~$t$.

We assume at first, then, that $D_\bdv$ is a nonstandard boundary
divisor. Since all the divisors of given $t$-divisor are equivalent under the
action of $\ZZ_2\wr \SL(2,\ZZ_t)$, (because the $t$-divisor is the only
invariant of a boundary divisor of $\cA_t$: see \cite{FS}) it will be
enough to calculate the number of conditions imposed by one divisor of
each type. That is to say, we only need consider boundary components
in $\cA_t^*$.

In view of this we may take $\bdv=(0,0,r,1)$ for some $r\vert t$ with
$1<r<t$. We write $(0,0,0,1)=\vz$ (for consistency with \cite{HKW2})
and we put $h=t/r$. Since we want to work with $\Gtbl$ rather than $\Gtblti$
(so as to use fractional linear transformations) we must consider the lines
$\QQ\,\bdv R_t=\QQ\,\bdv'$, where $\bdv'=(0,0,1,h)$, and
$\QQ\,\vz R_t=\QQ\,\vz$.

Note that $\bdv'Q_r=\vz$, where
\begin{equation*}
Q_r=\begin{pmatrix}1&1&0&0\\ h-1&h&0&0\\ 0&0&h&1-h\\
0&0&-1&1\end{pmatrix}
\in\Sp(4,\ZZ).
\end{equation*}
\begin{proposition}\label{FJnonstandard}
If $\bdv$ has $t$-divisor $r\neq t$, $1$, and $F\in\gothS_k^*(\Gtbl)$
is a cusp form of weight~$k$, then there are coordinates $\tau_i^\bdv$
such that $F$ has a Fourier expansion near $D_\bdv$ as
\begin{equation*}
F=\sum_{w\geq 0}\theta_w^\bdv(\tau_1^\bdv,\tau_2^\bdv)
\exp{2\pi iw\tau_3^\bdv/rt}.
\end{equation*}
\end{proposition}

\begin{proof}
As usual (cf.~\cite{HKW2}) we write $\cP_\bdv'$ for the stabiliser of
$\bdv'$ in $\Sp(4,\RR)$, so $\cP_\bdv'=Q_r^{-1}\cP_\vz Q_r$. We take
$P_\bdv'=\cP_\bdv'\cap\Gtbl$: this group determines the structure of
$\AtblV$ near~$D_\bdv$. It is shown in~\cite[Proposition~I.3.87]{HKW2}
that $\cP_\vz$ is generated by $g_1(\gamma)$ for
$\gamma\in\SL(2,\RR)$, $g_2=\zeta$, $g_3(m,n)$ and $g_4(s)$ for $m$,
$n$, $s\in\RR$, where
\begin{equation*}
g_1(\gamma)=\begin{pmatrix}a&0&b&0\\0&1&0&0\\c&0&d&0\\0&0&0&1\end{pmatrix}
\quad\mbox{for}\quad
\gamma=\begin{pmatrix}a&b\\c&d\end{pmatrix}
\end{equation*}
and $g_3$ and $g_4$  are given by
\begin{equation*}
g_3(m,n)=\begin{pmatrix}1&0&0&n\\m&1&n&0\\0&1&0&-m\\0&0&0&1\end{pmatrix},
\qquad 
g_4(s)=\begin{pmatrix}1&0&0&0\\0&1&0&s\\0&0&1&0\\0&0&0&1\end{pmatrix}.
\end{equation*}
So $P_\bdv'$ includes the subgroup generated by all elements of the form
 $Q_r^{-1}g_i Q_r$
with $a,b,c,d,m,n,s\in\ZZ$ which lie in~$\Gtbl$. In
particular it includes the lattice $\{Q_r^{-1}g_4(rts) Q_r\mid
s\in\ZZ\}$. If we take $Z^\bdv=Q_r^{-1}(Z)$ for
$Z=\begin{pmatrix}\tau_1&\tau_2\\
\tau_2&\tau_3\end{pmatrix}$ then we obtain 
\begin{equation*}
Z^\bdv=
\begin{pmatrix}h^2\tau_1-2h\tau_2+\tau_3&-h(h-1)\tau_1+(2h-1)\tau_2-\tau_3\\
-h(h-1)\tau_1+(2h-1)\tau_2-\tau_3&(h-1)^2\tau_1-2(h-1)\tau_2+\tau_3
\end{pmatrix}.
\end{equation*}
One easily checks that 
\begin{equation*}
Q_r^{-1} g_4(rt)Q_r:Z^\bdv=\begin{pmatrix}\tau^\bdv_1&\tau^\bdv_2\\
\tau^\bdv_2&\tau^\bdv_3\end{pmatrix}\Mapsto 
\begin{pmatrix}\tau^\bdv_1&\tau^\bdv_2\\
\tau^\bdv_2&\tau^\bdv_3+rt\end{pmatrix}
\end{equation*}
and this proves the result.
\end{proof}

We define a subgroup $\Gamma(t,r)$ of $\SL(2,\ZZ)$ by
\begin{equation*}
\Gamma(t,r)=\left\{\begin{pmatrix}a&b\\ c&d\end{pmatrix}\left.\right\vert a\equiv
d\equiv 1\hbox{ mod }t,\ b\equiv 0 \hbox{ mod }t^2,\ c\equiv 0 \hbox{
mod }r\right\}.
\end{equation*}
\begin{lemma}\label{notorsion}
If $D_\bdv$ is nonstandard then $P_\bdv'$ is torsion-free.
\end{lemma}
\begin{proof}
The only torsion in $\Gtbl$ is $2$-torsion and a simple
calculation shows that if $\Eins_4\neq g\in\cP_\vz$ and
$g^2=\Eins_4$, then $Q_r^{-1}g Q_r\not\in \Gtbl$ for $r\neq 1$, $t$.
\end{proof}

\begin{proposition}\label{FJformnonstand}
If $D_\bdv$ is nonstandard and $F\in\gothS_k^*(\Gtbl)$ then
$\theta_w^\bdv(r\tau_1^\bdv,t\tau_2^\bdv)$ is a Jacobi form of
weight $k$ and index~$w$ for $\Gamma(t,r)$.
\end{proposition}
\begin{proof}
By direct calculation we find that $Q_r^{-1}g_1(\gamma) Q_r\in\Gtbl$
if $\gamma\in\Gamma(t,r)$ and $Q_r^{-1}g_3(rm,tn) Q_r\in\Gtbl$ for
$m$, $n\in\ZZ$. Using these two elements, another elementary
calculation verifies that the transformation laws for Jacobi
forms given in~\cite{EZ} are satisfied, since
\begin{equation*}
Q_r^{-1} g_3(rm,tn)Q_r:Z^\bdv\Mapsto
\begin{pmatrix}\tau_1^\bdv&\tau_2^\bdv+rm\tau_1^\bdv+tn\\
\tau_2^\bdv+rm\tau_1^\bdv+tn&\tau^\bdv_3+2rm\tau_2^\bdv+r^2m^2\tau_1^\bdv
\end{pmatrix}
\end{equation*}
and
\begin{equation*}
Q_r^{-1} g_1(\gamma)Q_r:Z^\bdv\Mapsto
\begin{pmatrix}\gamma(\tau_1^\bdv)&\tau_2^\bdv/(c\tau_1^\bdv+d)\\
\tau_2^\bdv/(c\tau_1^\bdv+d)&\tau^\bdv_3-c\tau_2^\bdv/(c\tau_1^\bdv+d).
\end{pmatrix}
\end{equation*}
\end{proof}

\begin{lemma}\label{indexofbdygroup}
The index of $\Gamma(t,r)$ in $\Gamma(1)$ is equal to $rt\phi_2(t)$ for
$r\neq 1$, $t$.
\end{lemma}
\begin{proof}
Consider the chain of groups
\[
\Gamma(1)=\SL(2,\ZZ)>\Gamma_0(t)>\Gamma_0(t)(r)>\Gamma(t,r)
\]
and the normal subgroup $\Gamma_1(t)\vartriangleleft\Gamma_0(t)$, where
\begin{eqnarray*}
\Gamma_0(t)&=&\left\{\gamma=\begin{pmatrix}a&b\\ c&d\end{pmatrix}
\in\SL(2,\ZZ)
\mid \begin{array}{l}a\equiv d\equiv 1\hbox{ mod } t,\\ b\equiv 0 
\hbox{ mod } t\end{array}\right\},\\
\Gamma_1(t)&=&\left\{\gamma=\begin{pmatrix}a&b\\ c&d\end{pmatrix}\in\SL(2,\ZZ)
\mid \begin{array}{l}a\equiv d\equiv 1\hbox{ mod }t,\\ 
b\equiv c\equiv 0\hbox{ mod }t\end{array}\right\},\\
\Gamma_0(t)(h)&=&\left\{\gamma=
\begin{pmatrix}a&b\\ c&d\end{pmatrix}\in\SL(2,\ZZ)
\mid \begin{array}{l}a\equiv d\equiv 1\hbox{ mod }t,\\
 b\equiv 0 \hbox{ mod } t,\ c\equiv 0 \hbox{ mod } h\end{array}\right\}.
\end{eqnarray*}
Thus $\Gamma_0(t)(r)$ is the kernel of reduction mod~$r$ in
$\Gamma_0(t)$. By Corollary~\ref{orderofSL2},
$[\Gamma(1):\Gamma_1(t)]=t\phi_2(t)$. By the exact sequence
\begin{equation*}
0\To\Gamma_1(t)\To\Gamma_0(t)\To \left\{\begin{pmatrix}1&0\\ \bar c&1
\end{pmatrix}\mid \bar c\in\ZZ_t\right\}\imic \ZZ_t\To 0
\end{equation*}
we have $[\Gamma_0(t):\Gamma_1(t)]=t$, and similarly
\begin{equation*}
0\To\Gamma_0(t)(r)\To\Gamma_0(t)\To \left\{\begin{pmatrix}1&0\\ \bar c&1
\end{pmatrix}\mid \bar c\in\ZZ_r\right\}\imic \ZZ_r\To 0
\end{equation*}
gives $[\Gamma_0(t):\Gamma_0(t)(r)]=r$.

To calculate $[\Gamma(t)(r):\Gamma(t,r)]$ we let $\Gamma_0(t)(r)$ act on
$\ZZ_t\cross\ZZ_{t^2}$ by multiplication on the right, i.e. by $\gamma:(x,y)
\to (ax+cy, bx+dy)$. The stabiliser of $(1,0)\in\ZZ_t\cross\ZZ_{t^2}$ is then
$\{\bar\gamma\in\Gamma_0(t)(r)\mid a\equiv 1 \mod t, b\equiv 0 \mod t^2\}$,
which is $\Gamma(t,r)$. On the other hand the orbit of
$(1,0)\in\ZZ_t\cross\ZZ_{t^2}$ is $\left\{(\bar a, \bar b)\in
\ZZ_t\cross\ZZ_{t^2}\mid \begin{pmatrix}a&b \\ c&d\end{pmatrix}
\in\Gamma_0(t)(r)\right\}$: that is, the set of possible first rows of a 
matrix in $\Gamma_0(t)(r)$ taken mod~$t$ in the first column and mod~$t^2$
in the second. This is evidently equal to $\{(1, tb')\mid b'\in\ZZ_t\}$, 
and hence of size~$t$. Thus $[\Gamma(t)(r):\Gamma(t,r)]=t$, which 
completes the proof.
\end{proof}

The standard case is only slightly different, but now there is
torsion.
\begin{proposition}\label{FJformstand}
If $D_\bdv$ is standard and $F\in\gothS_k^*(\Gtbl)$ then
$\theta_w^\bdv(r\tau_1^\bdv,t\tau_2^\bdv)$ is a Jacobi form of
weight $k$ and index~$w$ for a group $\Gamma'(t,r)$, which
contains $\Gamma(t,r)$ as a subgroup of index~$2$.
\end{proposition}
\begin{proof}
Although the standard boundary components are most obviously given by
$(0,0,0,1)$ for $r=t$ and $(0,0,1,0)$ for $r=1$, we choose to take
advantage of the calculations that we have already performed by
working instead with $(0,0,t,1)$ and
$(0,0,1,1)$. Lemma~\ref{FJformnonstand} is still true, but we also
have $Q_t^{-1}\zeta Q_t\in \Gtbl$ and
$Q_1^{-1}(-\zeta)Q_1\in\Gtbl$. These give rise to the stated extra invariance.
\end{proof}
\begin{lemma}\label{numberofFJ}
The dimension of the space $J_{3k,w}\big(\Gamma'(t,r)\big)$ of Jacobi
forms of weight~$3k$ and index~$w$ for $\Gamma'(t,r)$ is given as a
polynomial in~$k$ and~$w$ by
\begin{equation*}
\dim J_{3k,w}\big(\Gamma'(t,r)\big)=
\delta rt\nu\left(\sfrac{kw}{2}+\sfrac{w^2}{6}\right)+\hbox{\rm linear terms}
\end{equation*}
where $\delta=\hf$ if $r=1$ or $r=t$ and $\delta=1$ otherwise.
\end{lemma}
\begin{proof}
By \cite[Theorem 3.4]{EZ} we have
\begin{equation}\label{Jsum}
\dim J_{3k,w}\big(\Gamma'(t,r)\big)\le 
\sum_{i=0}^{2w}\dim\gothS_{3k+i}\big(\Gamma'(t,r)\big).
\end{equation}
Since $\Gamma'(t,r)$ is torsion-free, the corresponding modular curve has
genus $1+\frac{\mu(t,r)}{12}-\frac{\nu(t,r)}{2}$, where $\mu(t,r)$ is the
index of $\Gamma'(t,r)$ in $\PSL(2,\ZZ)$ and $\nu(t,r)$ is the number of
cusps (see~\cite[Proposition 1.40]{Shi}). Hence by~\cite[Theorem 2.23]{Shi}
the space of modular forms satisfies
\begin{eqnarray}\label{dimofmodfms}
\dim\gothS_k\big(\Gamma'(t,r)\big)&=&k\left(\sfrac{\mu(t,r)}{12}
-\sfrac{\nu(t,r)}{2}\right)+\sfrac{k}{2}\nu(t,r)+O(1)\nonumber\\
&=&\frac{k\mu(t,r)}{12}+O(1)
\end{eqnarray}
as a polynomial in~$k$. By Lemma~\ref{indexofbdygroup} we have
$\mu(t,r)=\hf rt\phi_2(t)=rt\nu$ for the nonstandard cases,
$\mu(t,1)=\hf t\nu$ and $\mu(t,t)=\hf t^2\nu$. Now the result follows from
equations~\eqref{dimofmodfms} and~\eqref{Jsum}.
\end{proof}
If $F\in\gothS_{3k}^*(\Gtbl)$ then $F.(d\tau_1\wedge d\tau_2\wedge
d\tau_3)^{\tens k}$ extends over the component $D_\bdv$ if and only if
$\theta_w^\bdv=0$ for all $w<k$: see~\cite[Chapter IV, Theorem
1]{AMRT}. Hence the obstruction $\Omega_\bdv$ coming from the boundary
component $D_\bdv$ is
\begin{equation}\label{Omegav}
\Omega_\bdv=\sum_{w=0}^{k-1}\dim J_{3k,w}\left(\Gamma'(t,r)\right)
\end{equation}
where $\Gamma'(t,r)=\Gamma(t,r)$ if $D_\bdv$ is nonstandard.

By Corollary~\ref{numberofbdycpts} the total obstruction from the boundary is
\begin{equation*}
\Omega_\infty=\sum_{r|t}\#\Bar\Phi(h)\#\Bar\Phi(r)\sum_{w=0}^{k-1}\dim
J_{3k,w}\big(\Gamma'(t,r)\big),
\end{equation*}
and we may assume that $k$ is even. 
\begin{corollary}\label{bdyobstr}
The obstruction coming from the boundary is
\begin{equation*}
\Omega_\infty\le \bigg(\sum_{r\vert t}\delta rt\nu\#\Bar\Phi(h)
\#\Bar\Phi(r)\bigg)
\sfrac{11}{36}k^3+O(k^2).
\end{equation*}
\end{corollary}
\begin{proof}
Summing the expression in Lemma~\ref{numberofFJ} for $0\le w <k$, as
required by equation~\eqref{Omegav} gives the coefficient of
$\frac{11}{36}$ and the rest comes directly from
Lemma~\ref{numberofFJ} and Corollary~\ref{numberofbdycpts}.
\end{proof}

\section{Intersection numbers}

We need to know the degrees of the normal bundles of the curves that
generate $\Pic H_1$ and $\Pic H_2$. For this we first need to describe
the surfaces $H_1$ and $H_2$. The statements and the proofs are very
similar to the corresponding results for the case of $\Aplev$, given
in~\cite{HKW2} and~\cite{HS}. Therefore we simply refer to those sources
for proofs, pointing out such differences as there are.
\begin{proposition}\label{H1starisproduct}
$H_1$ is isomorphic to $X(t)\cross X(t)$. 
\end{proposition}
\begin{proof}
Identical to \cite[I.5.53]{HKW2}.
\end{proof}
\begin{proposition}
$H_2$ is the minimal resolution of a surface $\bar H_2$ which is given by
two $\SL(2,\ZZ_2)$-covering maps
\begin{equation*}
X(2t)\cross X(2t)\To \bar H_2 \To X(t)\cross X(t).
\end{equation*}
The singularities that are resolved are $\nu^2$ ordinary double points, one
over each point $(\alpha,\beta)\in X(t)\cross X(t)$ for which $\alpha$ and
$\beta$ are cusps.
\end{proposition}
\begin{proof}
Similar to~\cite[Proposition I.5.55]{HKW2} and the discussion
before~\cite[Proposition 4.21]{HS}. $X(2)$ and $X(2p)$ are both replaced
by $X(2t)$ and $X(1)$ and $X(p)$ by $X(t)$. Since $t>3$ there are no elliptic
fixed points and hence no other singularities in this case.
\end{proof}

\begin{proposition}\label{bdycurves}
$H^\circ_1$ and $H^\circ_2$ meet the standard boundary components
$D_\bdv$ transversally in irreducible curves $C_\bdv\imic X^\circ(t)$
and $C'_\bdv\imic X^\circ(2t)$ respectively. $D_\bdv$ is isomorphic to
the (open) Kummer modular surface $K^\circ(t)$, $C_\bdv$ is the zero
section and $C'_\bdv$ is the $3$-section given by the $2$-torsion
points of the universal elliptic curve over~$X(t)$.
\end{proposition}
\begin{proof} 
This is essentially the same as~\cite[Proposition I.5.49]{HKW2},
slightly simpler in fact. We may work with $\bdv=(0,0,1,0)$ and copy
the proof for the central boundary component in $\Aplev$, replacing
$p$ by $t$ (again the fact that $p$ is prime is not used).
\end{proof}

We do not claim that the closure of $D_\bdv$ is the Kummer modular
surface $K(t)$. They are, however, isomorphic near $H_1$ and $H_2$. We
remark that $H_1$ and $H_2$ do not meet the nonstandard boundary
divisors, because of Lemma~\ref{notorsion}.

\begin{proposition}
$\AtblV$ is smooth near $H_1$ and $H_2$.
\end{proposition}
\begin{proof}
Certainly $\Atbl$ is smooth since the only torsion in $\Gtbl$ is
$2$-torsion fixing a divisor in $\HH_2$. There can in principle be
singularities at infinity, but such
singularities must lie on corank~$2$ boundary components not meeting
$H_1$ nor $H_2$ (again this follows from Lemma~\ref{notorsion}).
\end{proof}
\begin{corollary}\label{Hsdisjoint}
$H_1$ does not meet $H_2$.
\end{corollary}
\begin{proof} 
Since $\AtblV$ and the divisors $H_1$ and $H_2$ are smooth at the
relevant points, the intersection must either be empty or contain a
curve. However, the intersection also lies in the corank~$2$ boundary
components. These components consist entirely of rational curves, and
if $t>5$ then $H_1\imic X(t)\cross X(t)$ contains no rational
curves. Hence $H_1\cap H_2=\emptyset$.

With a little more work one can check that this is still true for
$t\le 5$, but we are in any case not concerned with that.
\end{proof}

\begin{proposition}\label{PicH1}
The Picard group $\Pic H_1$ is generated by
the classes of $\Sigma_1=\bar C_{0010}$ and $\Psi_1=\bar
C_{0001}$. The intersection numbers are $\Sigma_1^2=\Psi_1^2=0$,
$\Sigma_1.\Psi_1=1$ and $\Sigma_1.H_1=\Psi_1.H_1=-\mu/6$.
\end{proposition}

\begin{proof}
As in~\cite[Proposition 4.18]{HS} (but one has to use the alternative
indicated in the remark that follows).
\end{proof}

\begin{proposition}\label{PicH2}
The Picard group $\Pic H_2$ is generated by
the classes of $\Sigma_2$ and $\Psi_2$, which are the inverse images of
general fibres of the two projections in $X(t)\cross X(t)$,
and of the exceptional curves $R_{\alpha\beta}$ of the
resolution $H_2\to \bar H_2$.
The intersection numbers in $H_2$ are
$\Sigma_2^2=\Psi_2^2=\Sigma_2.R_{\alpha\beta}=\Psi_2.R_{\alpha\beta}=0$,
$R_{\alpha\beta}.R_{\alpha'\beta'}=
-2\delta_{\alpha\alpha'}\delta_{\beta\beta'}$ and $\Sigma_2.\Psi_2=6$.
In $\AtblV$ we have $\Sigma_2.H_2=\Psi_2.H_2=-\mu$
and $R_{\alpha\beta}.H_2=-4$.
\end{proposition}

\begin{proof}
The same as the proofs of~\cite[Proposition 4.21]{HS}
and~\cite[Lemma 4.24]{HS}. The curves $R'_{(a,b)}$
from~\cite{HS} arise from elliptic fixed points so they are absent
here.
\end{proof}
Notice that $\Sigma_2$ and $\Psi_2$ are also images of the general
fibres in $X(2t)\cross X(2t)$ and are themselves isomorphic to~$X(2t)$.

\section{Branch locus}

The closure of the branch locus of the map $\HH_2\to\Atbl$ is $H_1\cup
H_2$ and modular forms of weight $3k$ (for $k$ even) give rise to
$k$-fold differential forms with poles of order $k/2$ along $H_1$ and
$H_2$. We have to calculate the number of conditions imposed by these
poles.
\begin{proposition}\label{H1obstr}
The obstruction from $H_1$ to extending modular forms of weight $3k$ to
$k$-fold holomorphic differential forms is
\begin{equation*}
\Omega_1\le
\nu^2\Big(\hf-\sfrac{7t}{24}+t^2\big(\sfrac{1}{24}
+\sfrac{1}{864}\big)\Big)k^3
+O(k^2).
\end{equation*}
\end{proposition}
\begin{proof}
If $F$ is a modular form of weight $3k$ for $k$ even, vanishing to
sufficiently high order at infinity, and
$\omega=d\tau_1\wedge d\tau_2\wedge d\tau_3$, then $F\omega^{\tens k}$
determines a section of $kK+\frac{k}{2}H_1+\frac{k}{2}H_2$, where $K$
denotes the canonical sheaf of $\AtblV$. From
\begin{equation*}
0\To \cO(-H_1)\To \cO\To \cO_{H_1}\To 0
\end{equation*}
we get, for $0\le j< k/2$
\begin{eqnarray*}
0&\To& H^0\big(kK+(\sfrac{k}{2}-j-1)H_1+\sfrac{k}{2}H_2\big)
\To H^0\big(kK+(\sfrac{k}{2}-j)H_1+\sfrac{k}{2}H_2\big)\\
&\To& H^0\big(\big( kK+(\sfrac{k}{2}-j)H_1+\sfrac{k}{2}H_2\big)\vert_{H_1}\big)
\end{eqnarray*}
so 
\begin{eqnarray*}
h^0\big(kK+(\sfrac{k}{2}-j)H_1+\sfrac{k}{2}H_2\big)&\le&
h^0\big(kK+(\sfrac{k}{2}-j-1)H_1+\sfrac{k}{2}H_2\big)\\
&&\mbox{}+h^0\big(\big(kK+(\sfrac{k}{2}-j)H_1+\sfrac{k}{2}H_2\big)\vert_{H_1}\big).
\end{eqnarray*}
Note that, by Lemma~\ref{Hsdisjoint}, $H_2\vert_{H_1}=0$.
Therefore 
\begin{equation*}
h^0\big(kK+\sfrac{k}{2}H_2\big)\ge
h^0\big(kK+\sfrac{k}{2}H_1+\sfrac{k}{2}H_2\big) + \sum_{j=0}^{k/2-1}
h^0\big(\big(kK+(\sfrac{k}{2}-j)H_1\big)\vert_{H_1}\big),
\end{equation*}
so
\begin{equation}\label{Omega1}
\Omega_1\le \Sum_{j=0}^{k/2-1}h^0\big(\big(kK
+(\sfrac{k}{2}-j)H_1\big)\vert_{H_1}\big)
=\Sum_{j=0}^{k/2-1} h^0\big(kK_{H_1}-(\sfrac{k}{2}+j)H_1\vert_{H_1}\big).
\end{equation}
By Lemma~\ref{PicH1}, $K_{H_1}$ and $H_1\vert_{H_1}$ are both
multiples of $\Sigma_1+\Phi_1$, and any positive multiple of
$\Sigma_1+\Psi_1$ is ample. Suppose 
$H_1\vert_{H_1}=a_1(\Sigma_1+\Psi_1)$ and $K_{H_1}=b_1(\Sigma_1+\Psi_1)$.
Then
\begin{equation*}
-\frac{\mu}{6}=\Sigma_1.H_1=a\Sigma_1.(\Sigma_1+\Psi_1)=a_1
\end{equation*}
and 
\begin{equation*}
\frac{\mu}{6}-\nu=2g(\Sigma_1)-2=(K_{H_1}+\Sigma_1).\Sigma_1=
K_{H_1}.\Sigma_1=b_1
\end{equation*}
Hence, using equation~\eqref{Omega1}
\begin{eqnarray*}
\Omega_1&\le&
\Sum_{j=0}^{k/2-1}
h^0\big((\sfrac{k\mu}{6}-k\nu+\sfrac{k\mu}{12}+\sfrac{j\mu}{6})
(\Sigma_1+\Psi_1)\big)\\
&=&\Sum_{j=0}^{k/2-1}
h^0\big((\sfrac{kt\nu}{4}-k\nu+\sfrac{jt\nu}{6})
(\Sigma_1+\Psi_1)\big).
\end{eqnarray*}
Since $t\ge 7$ (we know from~\cite{Muk} that $\AtblV$ is rational for
$t\le 5$), we have
$\sfrac{kt\nu}{4}-k\nu+\sfrac{jt\nu}{6}-\sfrac{t\nu}{6}+\nu>0$ for
all~$j$ and hence $(\sfrac{kt\nu}{4}-k\nu+\sfrac{jt\nu}{6})
(\Sigma_1+\Psi_1)-K_{H_1}$ is ample. So by vanishing we have
\begin{eqnarray*}
\Omega_1&\le& \Sum_{j=0}^{k/2-1}\hf
(\sfrac{kt\nu}{4}-k\nu+\sfrac{jt\nu}{6})^2
(\Sigma_1+\Psi_1)^2+O(k^2)\\
&=&\Sum_{j=0}^{k/2-1}(\sfrac{kt\nu}{4}-k\nu+\sfrac{jt\nu}{6})^2+O(k^2)\\
&=&\nu^2\big(\hf-\sfrac{7t}{24}+t^2(\sfrac{1}{24}+\sfrac{1}{864})\big)k^3
+O(k^2).
\end{eqnarray*}
\end{proof}
Next we carry out the same calculation for $H_2$.
\begin{proposition}\label{H2obstr}
The obstruction from $H_2$ is
\begin{equation*}
\Omega_2\le
\nu^2\big((\sfrac{1}{2}+\sfrac{1}{72})t^2-(\sfrac{1}{4}+\sfrac{1}{24})t
-\sfrac{7}{3}+\sfrac{1}{24})k^3
+O(k^2).
\end{equation*}
\end{proposition}
\begin{proof}
By the same argument as above (equation~\eqref{Omega1}) the obstruction is
\begin{equation*}
\Omega_2\le \Sum_{j=0}^{k/2-1}
h^0\big(kK_{H_2}-(\sfrac{k}{2}+j)H_2\vert_{H_2}\big).
\end{equation*}
In this case $H_2\vert_{H_2}=a_2(\Sigma_2+\Psi_2)+c_2R$, where
$R=\sum_{\alpha,\beta}R_{\alpha\beta}$ is the sum of all the exceptional
curves of $H_2\to \bar H_2$, and $K_{H_2}=b_2(\Sigma_2+\Psi_2)+d_2R$.
Since $\Sigma_2\imic X(2t)$ we have by~\cite[1.6.4]{Shi}
\begin{equation*}
2g(\Sigma_2)-2=\sfrac{1}{3}(t-3)\nu(2t)=\mu-\sfrac{\nu}{2}.
\end{equation*}
Hence
\begin{equation*}
-\mu=\Sigma_2.H_2=a_2\Sigma_2^2+a_2\Sigma_2.\Psi_2+c_2\Sigma_2.R=6a_2
\end{equation*}
so $a_2=-\mu/6$, and
\begin{equation*}
-4\nu^2=R.H_2=a_2\Sigma_2.R+a_2\Psi_2.R+c_2 R^2=-2\nu^2c_2
\end{equation*}
so $c_2=2$. Therefore
\begin{equation*}
H_2\vert_{H_2}=-\sfrac{\mu}{6}(\Sigma_2+\Psi_2)+2R.
\end{equation*}
Similarly
\begin{equation*}
\mu-\sfrac{\nu}{2}=(K_{H_2}+\Sigma_2).\Sigma_2=6b_2
\end{equation*}
so $b_2=\mu/6-\nu/12$, and $0=R.K_{H_2}=d_2R^2$ so $d_2=0$. Hence
\begin{equation*}
K_{H_2}=\sfrac{1}{6}(\mu-\sfrac{\nu}{2})(\Sigma_2+\Psi_2).
\end{equation*}
Moreover $L_j=(k-1)K_{H_2}-(\sfrac{k}{2}+j)H_2\vert_{H_2}$ is ample,
as is easily checked using the Nakai criterion and the fact that the
cone of effective curves on $H_2$ is spanned by $R_{\alpha\beta}$ and by
the non-exceptional components of the fibres of the two maps
$H_2\to X(t)$. These components are
$\Sigma_\alpha\equiv\Sigma_2-\sum_\beta R_{\alpha\beta}$ and
$\Psi_\beta\equiv\Psi_2-\sum_\alpha R_{\alpha\beta}$, and it is simple to
check that $L_j^2$, $L_j.\Sigma_\alpha=L_j.\Psi_\beta$ and
$L_j.R_{\alpha\beta}$ are all positive for the relevant values of $j$, $k$
and~$t$. Therefore
\begin{eqnarray*}
\Omega_2&\le& \Sum_{j=0}^{k/2-1}\hf\big(kK_{H_2}
-(\sfrac{k}{2}+j)H_2\vert_{H_2}\big)^2\\
&=&\Sum_{j=0}^{k/2-1}\hf\big(\nu(\sfrac{kt}{4}-\sfrac{k}{12}+\sfrac{jt}{6})
(\Sigma_2+\Psi_2)+(k+2j)R\big)^2\\
&=&\nu^2k^3\big(t^2(\sfrac{3}{8}+\sfrac{1}{8}+\sfrac{1}{72})
-t(\sfrac{1}{4}+\sfrac{1}{24})+\sfrac{1}{24}-2-\sfrac{1}{3}\big)+O(k^2)
\end{eqnarray*}
since $(\Sigma_2+\Psi_2)^2=12$.
\end{proof}

\section{Final calculation}

In this section we assemble the results of the previous sections into
a proof of the main theorem.
\begin{theorem}\label{mainthm}
$\AtblV$ is of general type for $t$ odd and $t\geq 17$.
\end{theorem}
\begin{proof}
We put $n=3k$ in Theorem~\ref{cuspforms}, and use $\phi_2(t)=2\nu$ and
the fact that 
\begin{equation*}
\phi_4(t)=t^4\Prod_{p|t}(1-p^{-4})=t^2\phi_2(t)\Prod_{p|t}(1+p^{-2}).
\end{equation*}
This gives the expression
\begin{equation*}
\dim \gothS_n^*(\Gtbl)={\frac{k^3\nu^2}{320}}t^4\Prod_{p|t}(1+p^{-2})+O(k^2).
\end{equation*}
From Proposition~\ref{H1obstr} and Proposition~\ref{H2obstr} we have
\begin{eqnarray*}
\Omega_1&=&
k^3\nu^2\left(\sfrac{37}{864}t^2-\sfrac{7}{24}t+\sfrac{1}{2}\right)
+O(k^2),\\
\Omega_2&=& k^3\nu^2\left(\sfrac{37}{72}
t^2-\sfrac{7}{24}t-\sfrac{55}{24}\right)
+O(k^2)
\end{eqnarray*}
and from Corollary~\ref{bdyobstr} and Corollary~\ref{numberofbdycpts}
\begin{equation*}
\Omega_\infty=k^3\nu^2\Sum_{r|t}\sfrac{11}{36r}t^2
 \Prod_{p|(r,h)}(1-p^{-2})+O(k^2).
\end{equation*}
since $\phi_2(r)\phi_2(h)=t^2\Prod_{p|(r,h)}(1-p^{-2})$.

It follows that $\AtblV$ is of general type, for odd~$t$, provided
\begin{equation}\label{ineq}
\frac{1}{320}\Prod_{p\vert t}(1+p^{-2})t^4
-\frac{481}{864}t^2+\frac{7}{12}t+\frac{43}{24}
-\Sum_{r|t}\frac{11}{36r}t^2\Prod_{p|(r,h)}(1-p^{-2})>0.
\end{equation}
This is simple to check: since either $r=1$ or $r\ge 3$, and since the
sum of the divisors of $t$ is less than $t/2$, the last term can be
replaced by $-\frac{11}{36}t^2-\frac{11}{108}t^3$ and the $t$ and
constant terms, and the the $p^{-2}t^4$ term, can be discarded as they
are positive. The resulting expression is a quadratic in $t$ whose
larger root is less than~$40$, so we need only consider odd $t\le
39$. We deal with primes, products of two primes and prime powers
separately. In the case of primes, the expression on the left-hand
side of the inequality~\eqref{ineq} becomes
$\frac{1}{320}t^4-\frac{7433}{8640}t^2+\frac{5}{18}t+\frac{43}{24}$,
which is positive for $t\ge 17$. The expression in the case of $t=pq$
is positive if $t\ge21$. For $t=p^2$ we get an expression which is
negative for $t=9$ but positive for $t=25$, and for $t=p^3$ the
expression is positive.
\end{proof}

One can say something even for $t$ even, though not if $t$ is
a power of~$2$.
\begin{corollary}
$\AtblV$ is of general type unless $t=2^ab$ with $b$ odd and $b<17$.
\end{corollary}
\begin{proof}
${\cA_{nt}^{\mathrm{bil}}}$ covers $\Atbl$ for
any~$n$, and therefore ${\cA_{nt}^{\mathrm{bil*}}}$ is of general type
if $\AtblV$ is of general type. 
\end{proof}

\begin{tabular}{l}
G.K. Sankaran\\
Department of Mathematical Sciences\\
University of Bath \\
Bath BA2 7AY\\
England\\
\\
{\tt gks@maths.bath.ac.uk}
\end{tabular}

\end{document}